\documentclass[11pt,twoside]{article}

\usepackage{amsmath}
\usepackage{amssymb}
\usepackage{graphicx}
\usepackage{color}

\newtheorem{theo}{Theorem}[section]
\newtheorem{rema}[theo]{Remark}

\newtheorem{propo}[theo]{Proposition}

\makeatother
\font\ddpp=msbm10  at 11 truept
\def\R{\hbox{\ddpp R}}
\def\In{\hbox{\ddpp I}}

\def\C{\hbox{\ddpp C}}

\def\Z{\hbox{\ddpp Z}}

\def\N{\hbox{\ddpp N}}
\def\E{\hbox{\ddpp E}}
\def\LL{\hbox{\ddpp L}}
\def\1{\mathbf 1}

\def\semidir{\ltimes}

\def\qed{\ensuremath{\quad\Box\quad}}
\def\ov{\overline}

\newcommand{\beq}{\begin{equation}}
\newcommand{\eeq}{\end{equation}}

\newcommand{\ben}{\begin{enumerate}}
\newcommand{\een}{\end{enumerate}}

\newcommand{\bit}{\begin{itemize}}
\newcommand{\eit}{\end{itemize}}

\newcommand{\bean}{\begin{eqnarray*}}
\newcommand{\eean}{\end{eqnarray*}}
\newcommand{\benu}{\begin{enumerate}}
\newcommand{\eenu}{\end{enumerate}}
\newcommand{\eea}{\end{eqnarray}}
\newcommand{\bea}{\begin{eqnarray}}

\font\ddpp=msbm10 scaled \magstep 1 
\def\R{\hbox{\ddpp R}}    

\def\N{\hbox{\ddpp N}}    
\def\v{\noindent}

\newcommand{\eart}{\end{article}}
\newcommand{\edoc}{\end{document}}

\begin{document}


\title{An Invitation to Lorentzian Geometry}%

\author{Olaf M\"uller\footnote{Fakult\"at f\"ur Mathematik der Universit\"at Regensburg, Universit\"atsstra\ss e 31, 93053 Regensburg (Germany). http://homepages-nw.uni-regensburg.de/~muo63888, \ \ olaf.mueller@mathematik.uni-regensburg.de}  \ and Miguel S\'anchez\footnote{Miguel S\'anchez Caja, Departamento de Geometr\'ia y Topolog\'ia, Facultad de Ciencias, Universidad de Granada. Campus de Fuentenueva s/n. E-18071 Granada (Spain). http://gigda.ugr.es/sanchezm/index.html, sanchezm@ugr.es}}

\date{}%
\maketitle

\begin{abstract}
The intention of this article is to give a flavour of some global problems in General
Relativity. We cover a variety of topics, some of them related to the fundamental concept of {\em Cauchy hypersurfaces}: (1)
structure of globally hyperbolic spacetimes, (2)
the relativistic initial value problem, (3) constant mean curvature surfaces, (4) singularity theorems, (5) cosmic
censorship and Penrose inequality, (6) spinors and holonomy.
\end{abstract}

\small{\tiny {\em Keywords}: Global Lorentzian Geometry, Cauchy
hypersurface, global hyperbolicity,  Einstein equation, initial
value problem, CMC hypersurface, singularity theorems, ADM mass,
cosmic censorship hypotheses, Penrose inequality, spinors,
Lorentzian holonomy.

}

\section{Introduction} 

The ground-breaking discovery of the theory of relativity has
revealed that the physics of gravity can be described successfully
by a theory treating space and time on the same footing,
distinguished by the sign of an inner product. Nowadays, the
mathematical framework of General Relativity can be regarded as a
branch of  Geometry (Lorentzian Geometry), in a similar sense like
mathematics of Theoretical Mechanics are a branch of symplectic
geometry. Admittedly, the physical leading ideas behind the
geometric results are more subtle and less evident in General
Relativity than in Mechanics. But after getting some familiarity
with it, a new geometric world opens up, including unexpected new
solutions (and problems) in {\em Riemannian} Geometry. We want to
provide a brief overview of this wonderful world --- which might
be closer to the reader's field of research than (s)he may expect.
We focus on global problems, which are usually the most
interesting for mathematicians. We hope that this article will be
also of some interest for physicists, which, frequently, are very
familiar with local Differential Geometry, while  global problems
are neglected as of little importance for experimental purposes.
Nevertheless, global questions can provide the necessary framework
to the full theory and have implications in more practical issues --- 
prominent examples supporting this claim are the definition
of mass and energy or the Aharonov-Bohm effect.

\bigskip

\vspace*{-1.0ex}
\begin{figure}[h]

\centering
\includegraphics[scale=8]{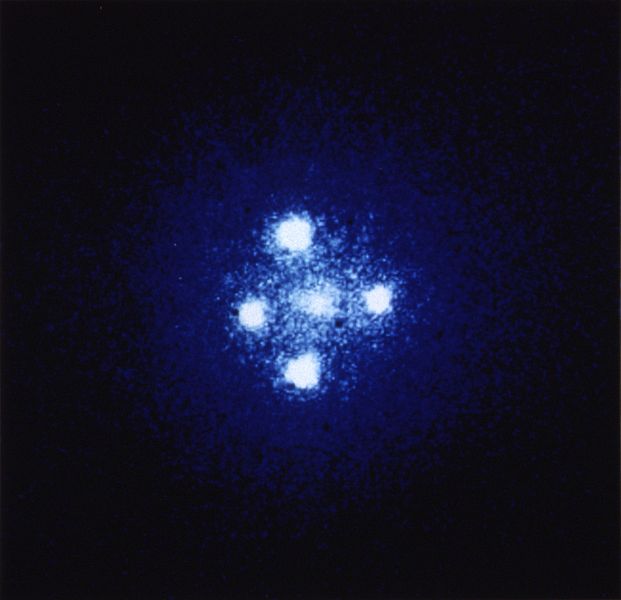}
\caption{{\footnotesize Einstein's cross, a gravitationally lensed region at 8 billion light years from earth, in the constellation Pegasus. A central galaxy focuses the light emanating from a quasar behind it producing a fourfold image of it. The images can be attributed to a single quasar via comparing characteristic spectra and synchronized relativistic jets. Gravitational lensing is a physical phenomenon closely related to the Lorentzian topic of Morse theory of causal curves \cite{vP} (see Section \ref{Geodesics}). Image taken by the European Space Agency's Faint Object Camera on board NASA's Hubble Space Telescope in 1990}}

\end{figure}

\vspace*{-0.5ex}


One can distinguish different
directions in Lorentzian Geometry: \ben \item
{\em Looking at Riemannian Geometry.} That is, one tries to adapt the
familiar Riemannian tools and results to the Lorentzian case,
as far as possible. This was the case at the beginning of General
Relativity, and is also the typical starting point for a standard
mathematician --- who has studied Riemannian Geometry but not
Lorentzian geometry. This is not as straightforward as it sounds, because
Lorentzian and Riemannian geometries, in spite of sharing common
roots, separate fast in both aims and methods.

\item {\em Developing specific Lorentzian tools}. Concepts such as
causality, boundaries, conformal extensions (Penrose diagrams),
asymptotic behaviors (spatial and null infinities) or black holes,
are specific to Lorentzian Geometry, without any
analog in the Riemannian case. Here, physical intuitions are  a
very important  guide, but we emphasize that these concepts have a
completely tidy mathematical definition.

\item {\em Feed back to Riemannian Geometry}. Sometimes, a problem
in Lorentzian Geometry admits a full reduction to a purely
Riemannian problem. This problem may be unexpected from a
Riemannian approach, but now it becomes  natural. The initial
value constraint equations for Einstein's equation, the positive
mass theorems (which yield the last step in the solution to the Yamabe
problem!) or the Penrose inequalities provide remarkable examples
of this situation. The reader will be able to appreciate that all
the main results in Section 7 are stated in a purely Riemannian
way, even though some motivations in this section, as well as
further developments in Section 8, show the power and beauty of
the bigger Lorentzian world. \een In what follows, after a
preliminary comparison between Lorentzian and Riemannian
Geometries, a short overview of six research areas in Lorentzian
Geometry, most of them motivated by General Relativity, is
provided. In our choice of problems, global hyperbolicity and
Cauchy hypersurfaces  play an important role. The reason is
twofold: on one hand, they play a central role in global problems,
on the other hand, they are a very intuitive bridge between
Riemannian and Lorentzian geometry (see Section 3). Additionally,
the reader will find a variety of further topics, including
hyperbolic equations, geodesics, CMC hypersurfaces, mass
inequalities, holonomy and spinors. Some parts of this paper
extend and update the earlier review \cite{S}.

\section{From Riemannian to Lorentzian geometry}
\label{s2}

In this section, we briefly recall the basics of
Lorentzian geometry and compare them to the Riemannian
situation (see \cite{BEE, CanSan, HE, MinSan, O, Wa} for further
details). Let, throughout this article, $M$ be an $n$-dimensional manifold, oriented if necessary. 
A {\em Lorentzian metric on $M$} is a symmetric bilinear form of
signature $(1,n-1)$ at every point of $M$, that is, the maximal
dimension of a negative definite subspace of $T_pM$ is $ 1 $ and
the maximal dimension of a positive definite subspace of $T_pM$ is
$ n-1 $. While on any manifold there is a Riemannian metric (using
paracompactness), the same is not true any more if one replaces
``Riemannian'' by ``Lorentzian'': actually, the existence of a
Lorentzian metric on a manifold $M$ is easily seen to be
equivalent to the existence of a one-dimensional subbundle of the
tangent bundle which implies the vanishing of the Euler class of
$\tau_M$, or equivalently, of the Euler characteristic of $M$
\cite[Theorem 2.19]{hB}. For noncompact manifolds, however, this
yields no obstruction: every non-compact manifold carries a
Lorentzian metric.

In General Relativity, information is allowed to travel not along arbitrary curves, but only along {\em causal} ones. The notion of causality is first defined on the tangent spaces via the sign appearing in the Lorentzian metric. Namely (following the convention in \cite{MinSan}),
a non-zero tangent vector $v\in T_pM$ is
{\em causal} if it is either {\em timelike}, i.e. $g(v,v)<0$, or
{\em lightlike}, i.e. $g(v,v)=0, v\neq 0$. A vector is called {\em
spacelike} if $g(v,v) >0$ (in particular, this convention means that $v=0$ is
non-spacelike and also non-causal). Both lightlike vectors and
the 0 vector are called {\em null vectors}.

Let $J_g$ be the set of all the causal vectors and $I_g$ the set
of timelike vectors. Then it is easy to see that each of $I_g \cap T_pM$ and $J_g \cap
T_pM$ have two connected components. A connected Lorentzian
manifold $(M,g)$ is called {\em time-orientable} if $J_g$ is
disconnected; in this case, $J_g$ has exactly two connected
components (a fact that can be seen easily by general
connectedness arguments using lifts of curves into open fibers). A
 {\em time-orientation} is a choice of one component $J_g^+$,
which is called the {\em causal future} then (and its $J_g$-complement
$J_g^- := - J_g^+$ is called {\em causal past}). Analogous
considerations are valid for $I_g$. By a {\em spacetime},  we mean
a (connected) time-oriented Lorentzian manifold $(M,g)$ (the
choice of time-orientation is not denoted explicitly).

The notion of ``future'' is then transferred from $TM$ to $M$ by
curves: A (continuous) piecewise $C^1$ curve $c: I \subset \R
\rightarrow M$ is called {\em future-directed causal} (resp. {\em
future-directed timelike future}) iff $c'(t) \in J^+_g$ (resp.
$c'(t) \in I^+_g$) for all $t$ in each closed subinterval
$I_j\subset I$ where $c$ is $C^1$ --- and correspondingly for the
past. Then, we define the {\em chronological future} of $p$ as:
$$I_g^{+} (p) := \{ q \in M \vert \ \exists \  \hbox{ future-directed} \
  {\rm timelike \ } {\rm curve\ } c \ {\rm from \ } p \ {\rm to \
} q \ \}. $$ Analogously, the {\em causal future}
 $J_g^{+} (p)$ is defined replacing ``timelike'' by ``causal'' in
 previous definition, and adding by convention  $p\in J_g^{+}
 (p)$. There are natural dual 'past' notions, and the subscript $g$
 is removed when there is no possibility of confusion.
Trivially,
 $I^{\pm} (p) $ are always open,
but simple examples show that, in general, $J^{\pm} (p)$ are
neither open nor closed --- however, $\ov{I^{\pm} (p)} =
\ov{J^{\pm} (p)}$ always holds. A state of a classical relativistic system at a point $p$ can then only depend on events in the causal pastof $p$, and the state at $p$ in turn can only influence the physics in the causal future $J^+(p)$ of $p$.
The distinction of curves by their causal character
corresponds to different elements of physical reality: massive
bodies (e.g. observers) are supposed to travel along timelike
curves, massless particles along lightlike curves, whereas spacelike
curves do not admit a direct physical interpretation.

Many constructions in semi-Riemannian geometry are independent of
the signature of the metric. First of all, any (non-degenerate) metric
determines always a unique metric torsion-free connection
(Levi-Civita connection). The Riemannian curvature tensor is
defined in exactly the same way, and  geodesics are defined as
$\nabla$-autoparallel curves; in particular, local convexity (in
the sense of existence of a geodesically convex neighborhood
around every point) is ensured. In Lorentzian signature, timelike
geodesics are local maxima of an appropriately defined length
functional (see below) on causal curves, and this property can be
extended to lightlike geodesics, even though some  relevant
subtleties arise
(see, for example, 
\cite[Section 2]{MiSa}). In the framework of relativistic physics,
the ``observer'' corresponding to a future-directed timelike curve
$c$,  is considered to be falling freely when $c$ is a geodesic.
Future-directed lightlike geodesics are regarded as ``(trajectories
of) light rays''. No good interpretation holds for spacelike
curves, even if they are geodesics, except for very special
classes of spacetimes.

\newpage

\vspace*{-1.0ex}
\begin{figure}[h!]

\centering
\includegraphics[scale=0.3]{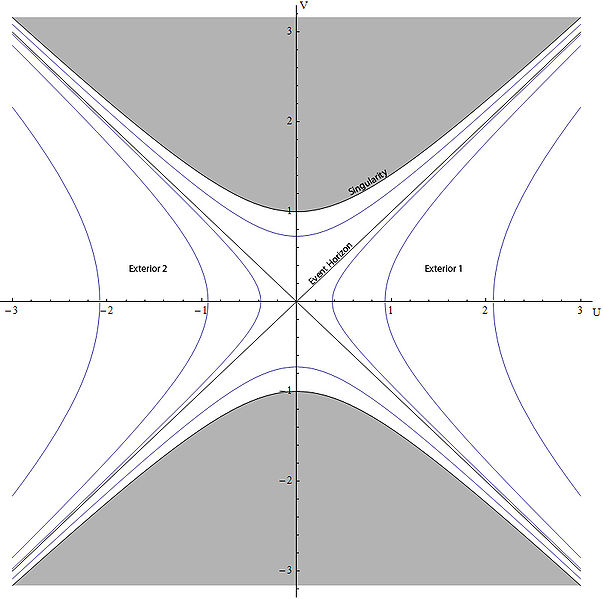}
\caption{{\footnotesize An instructive example: Schwarzschild-Kruskal metric. The two-dimensional Kruskal metric $h$ is conformally equivalent to the standard Minkowski metric on the non-shaded region $\{ (U,V) \in \mathbb{R}^2 \vert U^2 - V^2 > -1 \} $ (in particular, at every point, the future causal cone $J_g^+ \cap T_pM$ is just the upper quadrant between the translates of the diagonals). The conformal factor is $ \frac{32 M^3}{r} e^{-r/2M}$, and the $4$-dimensional physical spacetime is the warped product $h \times_{r^2} g_{\mathbb{S}^2}$ where $r$ is the unique positive number such that $U^2 - V^2 = (1 - \frac{r}{2M}) e^{r/2M}$. Lines through the origin are level sets of the coordinate $t$ which is a time function for the exterior of the event horizon, whereas the displayed hyperbola correspond to constant $r$ coordinate and are Killing orbits, timelike outside and spacelike inside the event horizon. The Schwarzschild part is exactly the future of the line $U = -V$; its exterior part describes successfully the gravitational field of an approximately spherically symmetric mass, e.g. of the sun. All future curves in the interior part have finite length (Image taken from http://commons.wikimedia.org/wiki/File:KruskalKoords.jpg licensed under the Creative Commons Attribution 3.0 Unported license by Author AllenMcC.).}}

\end{figure}

\vspace*{-0.5ex}

\bigskip

Many of the Riemannian constructions carry over to the
Lorentzian case {\em mutatis mutandis}. However, it is worthwhile
to give a brief overview of some of the most important differences
between Riemannian and Lorentzian geometry.


\begin{enumerate}
\item {\bf For any given manifold $M$, the set of all the
Lorentzian metrics on $M$ is not convex.} Recall that the set of
Riemannian metrics on $M$ is a (non-empty) convex cone. But this
does not hold for Lorentzian metrics, even when there is no
topological obstruction for their existence: one can check it
trivially in dimension 2 (if $g$ is Lorentzian then $-g$ is
Lorentzian too), and that counterexample can easily be extended to
higher dimensions. This fact makes several constructions of
interpolation between different metrics significantly more
complicated, e.g. in the theory of spinors.

\item{\bf Sectional curvature is defined
only for ``non-lightlike'' planes.} The reason is based on the fact
that when the restriction of the metric $g$  to some tangent
plane $\pi\subset T_pM$ becomes degenerate, then the denominator
in the definition of sectional curvature is 0. This fact also
has the consequence that, if the sectional curvatures of the non-degenerate
planes at  $p$ are not constant, then the sectional curvatures at
$p$ will reach values on all $\R$ (recall that  one would  divide
by arbitrarily small positive and negative quantities, depending
on the plane), see \cite{O}. Also due to elementary algebraic
reasons, a bound of the type Ric$(v,v)>cg(v,v)$ for all $v\in
T_pM\setminus 0$ and some $c\in \R$ cannot hold, and inequalities
for curvatures must be understood in the sense of Andersson and
Howard, \cite{AH}. In particular, the so-called {\em energy
conditions} (see Section \ref{s4}) will play the role of classical
curvature bounds in some Lorentzian results.

\item {\bf There is no good embedding of the category of
Lorentzian manifolds into metric measure spaces.} Recall that, in
order to understand collapse processes and prove finiteness
results like the ones due to Cheeger-Gromov, one needs to embed
the manifolds and their possible limits into some space of metric
spaces, independent of further specializations like Alexandrov
spaces, CAT(0) spaces, etc. There is a well-known natural
embedding of the category of Riemannian manifolds and convex open
isometrical embeddings to the category of metric measure spaces
which commutes with the forgetful functor to the category of
topological spaces (just by taking the geodesic distance and the
associated volume form). That is to say, from a Riemannian metric
one can construct in a natural way a metric space compatible with
the given topology. In contrast, there is no such construction in
the Lorentzian case (which is easily seen by the existence of
boosts around $p$ in Lorentz-Minkowski space $\LL^n:= \R^{1,n-1}$
mapping an arbitrary point in $J^+(p) \setminus I^+ (p)$ into an
arbitrarily small neighborhood of $p$). Consequently, all the
attempts of repeating the above compactness and finiteness results
are doomed to failure in the Lorentzian regime. The same holds for
most of the theory of isoperimetic inequalities. See, however, the
 the notion of ``Lorentzian distance'' below.

 \item{\bf The isotropy group of a point may be non-compact}. If a group $G$
acts faithfully and isometrically on a Riemannian manifold $(M,g)$,
then the  isotropy group of any $p \in M$, being a closed subgroup
of $O(n)$, must be compact. This does not hold by any means in the
Lorentzian case, as the Lorentz group $O_1(n)\equiv O(1,n-1)$
(which is the isotropy group of any $p\in \LL^n$) is not compact.
This represents an additional difficulty in finding invariant
quantities (which are usually found as an average w.r.t. an
associated Haar measure).

\item{\bf No analog to Hopf-Rinow holds}. Geodesic dynamics change
drastically: as there is no metric space associated to a given
Lorentzian metric,  none of the assertions of classical Hopf-Rinow
theorem hold in Lorentzian geometry. Geodesic completeness neither
implies b.a.-completeness (i.e., the completeness of curves with
bounded acceleration) as in the Riemannian case. Neither
compactness nor homogeneity of a Lorentzian manifold implies its
geodesic completeness (however, by an argument due to Marsden
\cite{Marsden}, one knows that both properties together do). A
counterexample for homogeneity is just a half-plane $H := \{ (x_0
, x_1) \in \LL^2 \vert x_0 > x_1  \}$ (the isometry group include
the actions of translations by $\R (\partial_0 +
\partial_1)$ plus the so-called boosts, i.e., the connected part of the identity of $O(1,1)$).
For compactness, consider the projection of the metric $g:= -2 du
dv + (cos^4 (v) - 1) du ^2)$ on the two-dimensional torus
$\R^2/\Z^2$. Putting $u:= x+t, v:= x-t$
the $g$-geodesic $t \mapsto (1/t - t, {\rm arctan} (t))$ is then
incomplete. Geodesics in Lorentzian manifolds can behave in a
strange way for the Riemannian intuition.  For example, there are
inextensible lightlike geodesics $c: I \rightarrow M$ that are
closed (in the sense that $c(I) = c(J)$ for a compact subinterval
$J$ of $I$) but they are non-periodic  (they appear in the
quotient of $H$ above by a discrete subgroup of boosts generated
e.g. by $ (u,v) \mapsto (2u, v/2)$ where $u,v$ are some natural
lightlike coordinates). Moreover, conjugate points along a
spacelike geodesic may accumulate \cite{Helfer, PicTausk}.

 \item {\bf Distinction between irreducibility and indecomposability for isometric actions}, as
 a Lorentzian vector space is not the sum of a (degenerate) subspace and its orthogonal. Recall that for an isometric Riemannian action $\rho$ on $T_pM$,
the orthogonal  of an invariant subspace $A$ of $\rho$, is
invariant itself and complements $A$. The latter ceases to be true
in the Lorentzian case if the metric restricted to $A$ is
degenerate. So, the action of a isometry group may be reducible to
$A$, but the vector space $T_pM$ may be not decomposable as sum of
irreducible parts. This applies to holonomy representations, and
made the discovery of the Lorentzian Berger list (the
classification of all possible holonomy groups) considerably more
difficult than the one of its Riemannian predecessor. This
landmark was completed in 2005, cf. Section \ref{Spinors}.

\end{enumerate}

At this point, Riemannian geometers should not feel
scared off by the above list of differences as they are balanced by a row of nice features listed below.
First of all, the Laplacian can be formally defined in the
Lorentzian setting as in the Riemannian one, but now it is an
hyperbolic operator (d'Alembertian) and, even more, {\em
Lorentzian Geometry is a natural framework to study hyperbolic
equations, as Riemannian Geometry is a natural setting for
elliptic ones}. This claim is supported by several geometric tools
available in the Lorentzian setting but not in the
Riemannian one. Let us point out some of these genuinely
Lorentzian tools.

\begin{enumerate}
\item {\bf Causality allows to visualize the conformal structure of  spacetimes.} In fact, the datum  $J_g^{\pm}$ which defines the causal structure is conformally invariant and, conversely, two Lorentzian metrics $g, g'$
are (pointwise) conformal (i.e., $g'= \Omega g$, $\Omega^2>0$)  if they have equal
causal cones. One can define the {\em chronological} $\ll$  and  {\em causal relations} $\leq$, namely
$p\ll q$ iff $p\in I^-(q)$, $p\leq q$ iff $p\in J^-(q)$.  Locally, either
of these relations characterizes the conformal structure.
So, {\em Causality} can be identified to conformal geometry  in
Lorentzian signature (nevertheless, a subtler  modification of the notion of Causality has
been recently introduced by Garc\'{\i}a-Parrado and Senovilla
\cite{GpSe}, see also \cite{GpSa}).
Remarkably, one can associate a (conformally invariant) {\em causal boundary} to every sufficiently well behaved spacetime, this allows to describe the possible asymptotics of timelike curves in a subtle way, ordering them by inclusion of their pasts \cite{FHS}.

\item {\bf Higher compatibility with conformal structures of
lightlike pregeodesics.} It is easy to  see that there are
positive functions of Schwartz class on Euclidean $\R^2$ which do
not leave, when used as a conformal factor, any pregeodesic
invariant; here, a pregeodesic means a geodesic up to a
reparametrization. In contrast, on a Lorentzian manifold,
conformal factors leave always some pregeodesics invariant, namely
{\em all the lightlike pregeodesics}. Moreover, conformal changes
even preserve conjugate points on them plus their multiplicities,
\cite{MinSan}.

\item {\bf Completeness and singularities.} Apart  from the
differences noted above, there are also new notions of
completeness in Lorentzian geometry. Traditionally, one
distinguishes (logically independent) weaker notions than geodesic
completeness: spatial completeness, lightlike completeness and
timelike completeness, depending on the causal character of the
geodesics in question. Such a subdivision of completeness does not
exist in Riemannian geometry and, even more, two stronger notions
appear in the Lorentzian setting: the above mentioned {\em
b.a.-completeness} and Schmidt's {\em b-completeness} \cite{Sc}
(any of them coincide with usual geodesic completeness in
Riemannian geometry). As a remarkable difference with the
Riemannian case, both, completeness and incompleteness are $C^r$
unstable, for every $r \in \N$, even in the case of Lorentzian metrics on a torus
\cite{RSjmp2}, but some results on stability can be still obtained
in the Lorentzian case, as $C^1$-fine stability in the globally hyperbolic case, see \cite{BEE, CanSan}.
The interplay between these particularities, causality and some physical
interpretations, have the effect that, in Lorentzian Geometry, {\em
singularity theorems} (which ensure incompleteness rather than
completeness) play an important role, see below.

\item {\bf Reverse triangle inequality and Lorentzian ``distance''}.
As a construction related to Causality but not conformally
invariant, one can define the {\em time separation} $d(p,q)$
between two points $p,q\in M$ of a spacetime $(M,g)$ as the {\em
supremum} of the lengths of the future-directed causal curves
starting at $p$ and ending at $q$; this supremum may be infinity,
and it is regarded as equal to 0 if $p\not\leq q$ (in fact
$d(p,q)=0$ iff $p\not\ll q$). The corresponding function
$d:M\times M\in [0,\infty]$ is commonly called the {\em Lorentzian
distance}, as it satisfies a {\em reverse} triangle inequality
(due to the existence of a reverse triangle
inequality for causal vectors in the same cone in any $T_pM$) with some
similarity to the Riemannian case. However, it presents also
some big differences with the Riemannian case; for example,
$d(p,q)>0$ implies either $d(q,p)=0$ or there exists a closed
timelike curve through $p$ and $q$ (and, so,
$d(p,p)=d(q,q)=\infty$). As a clear connection with Causality,
$d(p,q)>0$ iff $p\ll q$ and, under a mild Causality condition on
the spacetime (strong causality, i.e., absence of ``almost
closed''  causal curves) the Lorentzian metric $g$ can be
reconstructed from $d$, as in the Riemannian case. However, the
interplay of $d$ with Causality makes it a genuinely Lorentzian
element. For example, a spacetime is globally hyperbolic (see
definition below) iff the Lorentzian distance is finite and
continuous for all the metrics  in the conformal class of $g$ (see
for example \cite{BEE, MinSan}).

\end{enumerate}


\section{Causality and global hyperbolicity}

\v Good conditions on the causality of a spacetime may yield some
connections between Riemannian and Lorentzian manifolds and links
between hyperbolic and elliptic equations. A key notion is {\em
global hyperbolicity} which is to be developed here and which will play a
role in the spirit of {\em completeness} for Riemannian manifolds.

\bigskip

 A spacetime $(M,g)$ is
called {\em globally hyperbolic} iff it is
causal\footnote{Typically,  an a priori stronger hypothesis that
causality is used to define global hyperbolicity, namely, {\em
strong} causality. In \cite{BeSa4} it has been shown that both
notions agree. For simplicity, we renounce giving details here and
simply use the more recent definition of global hyperbolicity
instead.} (i.e. $p \notin J^+ (p)$ for all $p \in M$) and
diamond-compact. Here, ``diamond-compact'' means $J^+ (p) \cap J^-
(q)$ compact for all $p,q \in M$. The condition of causality
corresponds to the existence of global solutions of natural linear
differential operators for initial values on maximal achronal
hypersurfaces while the condition of diamond-compactness
corresponds to their uniqueness: consider the examples of a flat
Lorentzian torus for non-existence and of a flat vertical strip in
Minkowski space for non-uniqueness. In terms of physics,
diamond-compactness corresponds to predictability of nature or
Laplace's demon principle (which has been of some influence at least in
classical physics), whereas causality corresponds to the exclusion of
the possibility of time machines\footnote{Time machines are in
contradiction to the unspoken fundamental assumption of the free
will of the experimentalist taken by the vast majority of
physicists, in the sense that any observer, in contrast to
physical nature around him, is assumed to be able to take
decisions like preparing a spin-up or a spin-down state in a
manner which is in principle unpredictable for others, compare the
discussion of Bell's inequality and the EPR paradox. Note that
without that assumption, time machines do not contradict any other
principles of physics --- with the possible exception of
predictability of nature, cf. the article of Krasnikov \cite{sK}
as well as its critical reception in \cite{jbM}.}. The link
between global hyperbolicity and Riemannian completeness comes
from the following result, which lies in the spirit of
Hopf-Rinow's (see for example \cite{BEE, O}):
\begin{propo} \label{tas}
In a globally hyperbolic spacetime $(M,g)$, the Lorentzian
distance $d$ is finite, continuous and satisfies the {\em Avez-Seifert
property}, i.e., for any pair of causally related distinct points
$p,q\in M$, ($p\leq q$) there exists a causal geodesic from $p$ to
$q$ with length equal to $d(p,q)$.
\end{propo}
There are many examples of globally hyperbolic manifolds:

\begin{enumerate}
\item{A  Lorentzian product $(\In \times N, - dt^2 + g_N)$ for an
interval $\In$, is globally hyperbolic iff  $g_N$ is a a complete
Riemannian metric on $N$; in particular, Lorentz-Minkowski spaces
are globally hyperbolic.} \item{{\em Narrower Cones Principle}: if
$(M,g)$ is globally hyperbolic and $ h$ is another Lorentzian
metric on $M$ with $J_h \subset J_g$, then $(M,h)$ is globally
hyperbolic as well; in particular, global hyperbolicity is
conformally invariant. This implies that, instead of the
Lorentzian products in item 1, one can equally consider {\em
Generalized Robertson-Walker spacetimes} (with complete $g_N$)
i.e, warped products $(\In \times N, - dt^2 + f(t) g_N)$ for some
positive function $f$.}

\item{If $(M,g)$ is globally hyperbolic and $A \subset M$ is a
causally convex open subset of $M$ in the sense that causal curves
cannot leave and then re-enter $A$, then $(A, g \vert_A)$ is
globally hyperbolic as well.} \item{Using convex neighborhoods, it
is easy to see that any point in any Lorentzian manifold has a
globally hyperbolic neighborhood.} \item{Global hyperbolicity of
Lorentzian metrics is a $C^0$-fine stable property in the space of
Lorentzian metrics which has been shown by the works of Geroch
\cite{rG70} (taking into account the progress made by Bernal and
S\'anchez in \cite{BeSa2}, see \cite{Sa13}); see also Lerner
\cite{deL} or the extended version on arxiv.org of Benavides and
Minguzzi \cite{BNM}.}
\end{enumerate}
This last point goes in the direction of Proposition \ref{tas},
i.e. the role of global hyperbolicity is related to Riemannian
completeness, as both properties are $C^0$ stable (but geodesic
completeness is not $C^r$-stable for any $r$ in the general Lorentzian case and only $C^1$-fine-stable for globally hyperbolic manifolds, as mentioned above).

Geroch \cite{rG70} showed in 1970 that a spacetime is globally
hyperbolic\footnote{He considered a different, but equivalent,
notion of global hyperbolicity, based on the compactness of the
space of causal curves connecting each two points, but this is not
specially relevant at this point.} if and only if it contains a
Cauchy hypersurface, that is, a set $\Sigma$ that is crossed
exactly once by any inextensible timelike curve (a posteriori,
$\Sigma$ must be then a topological hypersurface, see \cite{O}).
Moreover, he gave a construction of a continuous {\em Cauchy time
function} $t$, which means that $t$ increases strictly
monotonously along every causal future-directed curve, and is
surjective onto $\R$ along any inextensible causal future-directed
curve. His construction involved volumes of sets type $J^\pm (p)$
for a finite volume form. For a long time, it was not known, but
generally assumed, that $\Sigma$ could be taken as a smooth and
spacelike (non-degenerate) hypersurface and, even more, that  for
such a prescribed Cauchy hypersurface $\Sigma$  one could find
even a {\em Cauchy temporal function} vanishing on $\Sigma$. This
term is more special than the one before and denotes a {\em
smooth} function $t$ whose gradient satisfies $g({\rm grad} t,
{\rm grad} t ) < 0$ with ${\rm grad} t$ past-directed, plus the
surjectivity property above. (Note that not all smooth Cauchy time
functions are temporal: consider, on $\LL^{2} $, the function
$t(x_0, x_1) := (x_0 + x_1)^3$). Functions of this kind
automatically lead to metric splitings, that is, they imply that
the manifold is isometric to

\bea \label{esplit} (\R \times N, - f^2 \cdot dt^2 + g_t) \eea
where $f>0$ is a smooth function on $\R \times N$ and $g_t$ is a
smooth family of Riemannian metrics on the level sets of $t$, and
all level sets of $t$ are Cauchy. The interest in these questions
is obvious: on one hand, (smooth) spacelike Cauchy hypersurfaces
are the natural ones for initial data  (Einstein equation, Penrose
inequality..., see the next sections); on the other, the
orthogonal splitting is useful for many properties:  Morse Theory,
quantization, to find global coordinates, etc. Moreover, it also
leads to remarkable analytic results: not only adapted linear
symmetric hyperbolic systems (that is, those given by a
first-order differential operator on a vector bundle whose symbol
is positive-definite exactly on $I_g$), enjoy global existence and
uniqueness for arbitrary smooth initial values at a Cauchy hypersurface,
but the same is true for appropriate nonlinear equations like
Yang-Mills equations, as shown by Chrusciel and Shatah
\cite{CS00}. Physically, each temporal function $t$ determines in
a natural way not only  a one-parameter family of diffeomorphic
``physical spaces'' (the slices $t=$ constant), but also a Wick
rotation, obtained by inverting the sign on $\R \cdot {\rm grad}
t$ and leaving the orthogonal complement unchanged.

Sachs and Wu \cite[p. 1155]{SW} posed the existence of a {\em
smooth} Cauchy hypersurface in any globally hyperbolic spacetime
as a first open ``folk'' problem. Such a type of  problems cannot
be overlooked by physicists as minor questions of mathematical
rigor, as the requirements in the definition of global
hyperbolicity are plausible from the physical viewpoint, but the
assumption  of a splitting {\em a priori} of the spacetime as in
(\ref{esplit}) (the type of expression truly useful for several
physical purposes) would be totally unjustified. In a series of
papers published along 2003-2006, Bernal and S\'anchez
\cite{BeSa1, BeSa2, BeSa3} (see also \cite{Sa}) gave a full
solution by showing that a splitting as \eqref{esplit} can be
obtained and, then, any prescribed spacelike Cauchy hypersurface
can be chosen as the level $t=0$ of the splitting; their proof
used local convex coordinates patched together in a sophisticated
way. There has been, however, quite a few of interesting
developments since then. In 2011, M\"uller and S\'anchez
\cite{MuSa} solved the question of which Lorentzian manifolds are
isometrically embeddable in some $\LL^n$ (in the spirit of Nash's
theorem). With this aim, they proved, in particular, that any
globally hyperbolic spacetime admits a splitting as in
(\ref{esplit}) with an upper bounded function $f<1$ (this yielded
directly the isometric embeddability of all globally hyperbolic
spacetimes). Further properties on both,  the splitting (bounds
for curvature elements of the slices, flexibility) and the
isometric embedding in $\LL^N$ (closedness) were obtained then by
M\"uller \cite{Mu1, Mu2}. In 2012, Fathi and Siconolfi \cite{FS12}
proved the existence of a Cauchy temporal function in a class of
geometric spaces with a {\em cone structure} (which generalized
notably the class of globally hyperbolic spacetimes); their proof
involves tools from weak KAM theory. By taking into account the
progress along these decades (including old work by Seifert),
Chrusciel, Grant and Minguzzi \cite{CGM} have proved very recently
that, for some appropriate non-canonical choice of a volume form,
also the original functions defined by Geroch become $C^1$ (and
can be smoothed out further by local convolutions). The interplay
among these tools is an exciting matter of study \cite{Sa13}.

Summing up from a broad perspective, classical elementary results
as Proposition \ref{tas}, deeper structural results as splitting
\eqref{esplit}, and links with other parts of Differential
Geometry or Mathematical Physics (Morse theory, Geometric
Analysis, Cosmic Censorship, Wick rotation, Einstein and Yang
Mills equations...), show that Riemannian geometry is an
indispensable tool in the theory of globally hyperbolic manifolds,
but the study of the interplay between the two regimes has just
been initiated.

\section{Initial value problem} \label{s4}

Einstein's field equation can be written (in suitable units) as
\beq  \label{einstein} \mbox{Ric} - \frac{1}{2} S g = 8 \pi T.\eeq
Here, the geometric terms on the left hand side (Ricci tensor Ric,
scalar curvature $S$) are related to a symmetric 2-tensor on the
right hand side, the ``stress-energy'' $T$, which describes the
distribution of matter/energy.

More properly, we must emphasize that the unknown quantity is not
only the metric $g$ (with Ric and $S$): equations for $T$ must be
added to get a coupled system with (\ref{einstein}). Nevertheless,
we will assume for simplicity (in addition to dim($M$)=4, when
necessary) the following cases:

\bit \item Along this section, $T=0$ (vacuum), i.e.
(\ref{einstein}) becomes Ric$\equiv 0$.

\item In the next sections, solutions with $T$ non-determined
but satisfying only any of the (mild) ``energy conditions'' as:
(1) Weak: $T(v,v)\geq 0$ for any timelike $v$ (density energy is
nonnegative), (2) Dominant: $-T(v,\cdot )^\flat\equiv
-g^{ij}T_{jk}v^k $ is either future-directed causal or 0
for any future timelike $v$ (energy flow is causal), (3) Strong:
equivalent via Einstein equation to the timelike convergence
condition, Ric$(v,v)\geq 0$ for timelike $v$ (gravity, on average,
attracts). \eit
 The well-posedness of Einstein equation requires an input of
 initial data on a  3--manifold $\Sigma$ which permits to obtain a (``unique,
maximal'') spacetime (and eventually a $T$) such that $\Sigma$ is
embedded in $ M$ consistently with the initial data. The problem
is complicated:  a classical theorem such as Cauchy-Kovalevskaya's
is not applicable,  and, even more, in principle the system of
equations is not hyperbolic. Nevertheless, there exist a highly
non-trivial procedure --- based on the existence of {\em harmonic
coordinates} --- which allows one to find an equivalent
(quasi-linear, diagonal, second order) hyperbolic system. The
standard global result was obtained by Choquet-Bruhat and  Geroch
\cite{ChGe}:

 \begin{theo} \label{tGC} Let $(\Sigma,h)$ be a (connected)
Riemannian 3-manifold,  and $\sigma$ a symmetric two covariant
tensor on $\Sigma$ which satisfies the compatibility conditions of a second
fundamental form (Gauss and Codazzi eqns.)
 Then there exist a unique spacetime $(M,g)$ satisfying the following conditions:

(i) $\Sigma \hookrightarrow M$, consistently with $h, \sigma$
(i.e., $h=g|_\Sigma$ etc.)

(ii) Vacuum: Ric$\equiv 0$ (this can be extended  to any family
$T$ of natural divergence-free symmetric $2$-tensors, e.g. to
those coming from natural symmetric hyperbolic field theories).

(iii) $\Sigma$ is a Cauchy hypersurface of $(M,g)$.

(iv) Maximality: if $(M',g')$ satisfies (i)---(iii) then it is
isometric to an open subset of $(M,g)$.

\end{theo}
As suggested previously, the property (iii) becomes essential for
the well-posedness of the problem --- namely, the existence of a
solution spacetime can be proven because no timelike curve crosses
$\Sigma$ twice, and the uniqueness because all timelike curves
cross $\Sigma$ at least once.

\begin{rema} \label{rSCCC} {\em {\em  (SCCC)}. Even though the solution $(M,g)$ provided by Theorem \ref{tGC}
is maximal, it may be extensible as a spacetime, that is, $(M,g)$
may be isometric to an open proper subset of  another spacetime
$(\bar M, \bar g)$ --- even a vacuum one. In this case, $\Sigma$
cannot be a Cauchy hypersurface of the extension, and two
possibilities arise: (a) $(\bar M, \bar g)$ is not globally
hyperbolic or (b) the initial $\Sigma$ was not ``chosen
adequately'', as an input hypersurface for a whole physically
meaningful spacetime. Thus, an important question is how to characterize the (in)extendibility of $(M,g)$.

This question becomes extremely important in General Relativity
because physical intuition suggests that spacetime is
inextensible, but it suggests at the same time that it should be
predictable from initial data and, thus, globally hyperbolic.

The {\em Strong Cosmic Censorship Conjecture (SCCC)} asserts that,
 for generic physically reasonable data (including a ``good choice'' of $\Sigma$), $(M,g)$ is inextendible.
 Of course, a non-trivial problem of the conjecture, is to explain
 carefully what ``generic physically reasonable data'' means.
}\end{rema}

 A systematically studied problem is to characterize/classify the
solutions of (vacuum) Einstein equation. By using Theorem
\ref{tGC}, this is rather a purely Riemannian problem (roughly:
given data as, say, ($\Sigma, h$), classify the $\sigma$'s which
satisfy Gauss and Codazzi equations). There are two specially
important methods of solution (see \cite{Ba} for a detailed
exposition and \cite{CGP} for updated references): \bit \item {\em
Conformal}. Initial data are divided into two subsets: a subset of
{\em freely specified} conformal data (the conformal class of $h$,
a scalar field $\tau$, and a symmetric divergence free 2-tensor
$\tilde \sigma$), and a subset of {\em determined } data (a
function $\phi>0$, a vector field $W \in \chi (\Sigma )$), which
are derived from the free data by means of differential equations.
The interpretation and equations for these data vary with  two
types of conformal method (the method (A) or semi-decoupling,
whose origin goes back to Lichnerowicz \cite{Li}, and the method
(B) or conformally covariant). The problem is then to show if
there exists solutions for the equations of the determined data,
and classify them. \item {\em Gluing solutions}.    As a
difference with the conformal method, this is not a general one,
but it is very fruitful in relevant particular cases. Corvino and
Schoen \cite{Co, CoSc} glue any bounded region of an
asymptotically flat spacetime with the exterior region of a slice
of Kerr's --- this case becomes specially interesting as the ``no
hair theorems'' highlight Kerr spacetime at the final state of the
evolution of  a black hole. The useful gluing by Isenberg et al.
(\cite{Is2}, see also the initial data engineering  \cite{Is1})
constructs consistent initial data for Einstein equation from the
connected sum of previously obtained data  (for example,
construction of wormholes). \eit For the general conformal method,
the results depend on different criteria
--- topology of $\Sigma$, asymptotic behaviour, regularity
(analytic, smooth, H\"{o}lder class...), metric conformal class
(Yamabe type)... The most important one is the mean curvature $H$.
Essentially, when $H$ is constant almost all is known  (at least
if $T=0$); in fact, if $\Sigma$ is either compact without
boundary, or asymptotically flat or hyperbolic, it is  completely
determined which solutions exist (and they exist for all but
certain special cases). When $H$ is nearly constant there are many
results, but also many open questions; otherwise, there are very
few results.

\section{Constant mean curvature spacelike hypersurfaces}

The importance of (spacelike) hypersurfaces of constant mean
curvature (CMC) $H$ in a spacetime $(M,g)$, has been stressed
above in relation to the initial value problem, but they are also
important for other issues in General Relativity (see the survey
by Marsden and Tipler \cite{MT}). Here we will explain some
results on their existence and uniqueness. We will restrict to the
case when the hypersurface $\Sigma$ is spacelike, as it is easy to see that
a submanifold can extremize area only when its dimension is equal either to the
index or to the coindex of the metric (otherwise, the area may be
critical, but not extremal). Moreover, when $H=0$
the spacelike hypersurfaces are either maximal or neither maximize
or minimize area. Nevertheless, in both cases they are usually
called {\em maximal}, just like Riemannian minimal surfaces. About
the results on existence, we point out (see Gerhardt's book
\cite{Ge_book} for a detailed study):

\begin{itemize}
\item After the special case of Lorentz Minkowski (see below), the
first natural  problem to be considered is the construction of one
CMC hypersurface or, in general, a spacelike hypersurface $\Sigma$
with a {\em prescribed} mean curvature $H$, in a given spatially
compact globally hyperbolic spacetime. A relevant result of global
existence was due to Claus Gerhardt \cite{cG83} in 1983, under the
condition of existence of barriers. An upper resp. lower
$r$-barrier is a closed spacelike achronal hypersurface with mean
curvature $> r$ resp. $<r$. If there is an upper $r$-barrier
$\Sigma^+ $ and a lower $r$-barrier $\Sigma^- $, Gerhardt shows
that there is a CMC hypersurface of mean curvature $r$ in $I^+
(\Sigma^- ) \cap I^- (\Sigma^+)$. This is shown by solving the
Dirichlet problem for a given boundary curve by making some a
priori gradient estimates and, then, by applying a Leray-Schauder
fixed point theorem. To enhance the constructiveness in the last
part, Ecker and Huisken \cite{kEgH} used an evolutionary equation
in terms of the mean curvature flow starting at some Cauchy
hypersurface $\Sigma \subset I^+ (\Sigma^- ) \cap I^- (\Sigma^+)$
which provided a better control over the hypersurfaces --- for
example, it allows to fix all points of vanishing mean curvature
during the process. Even though they had to assume some additional
conditions (the timelike convergence condition and a more
technical structural monotonicity condition), Gerhardt
\cite{cG2000} refined Ecker and Huisken's flow method, showing
that such additional conditions were unnecessary. The improved
control allows to solve also related problems such as: (a) given a
prescribed point in $M$, construct a CMC hypersurface passing
through the point, or (b) given a compact Cauchy surface in $M$
find a compact CMC Cauchy surface with the same volume.


\item Removing spatial compactness, the next step is to consider
the existence of CMC hypersurfaces in asymptotically flat
spacetimes. Substantial contributions, specially in the maximal (eventually up to a compact subset) case, have been made by Bartnik \cite{rB1984}, Bartnik, Chrusciel
and Murchadha \cite{BCM}, and Ecker \cite{kE1993}. All three
articles assume an energy inequality and, moreover, a connection
between radial and time variables (called ``uniform interior
condition'' in the first and
 ``bounded interior geometry'' in the
second article). Roughly, both versions of the latter condition assert that the
deviations from Minkowski geometry propagate with subluminal
velocity --- so one can expect the condition to be true for massive
Klein-Gordon Theory, e.g. Again, the first article uses
Leray-Schauder's fix point theorem, while the other articles use
long-time convergence of the mean curvature flow, which provides a
better control on the surfaces.

\item The question of constructing a whole foliation by
CMC hypersurfaces was also studied. In the cited 1983 article by
Gerhardt  \cite{cG83}, he considered  a globally hyperbolic
spacetime with compact Cauchy hypersurfaces satisfying the
timelike convergence condition. Under these hypotheses, slices
with a CMC $H\neq 0$ are unique and, if there are two different
maximal slices, then both have to be totally geodesic and the region enclosed by them must be static ---
thus, the existence of two different hypersurfaces of CMC implies strong
obstructions. In fact, the mean curvature must increase
monotonously in foliations by CMC hypersurfaces.
 The timelike convergence
condition is replaced by the mere assumption of a lower bound to
the Ricci tensor on timelike vectors in a second article by
Gerhardt \cite{cG06}. This article also treats exclusively the
case of globally hyperbolic spatially compact spacetimes. Here,
the statement is the following: If for some sequence of Cauchy
surfaces $\Sigma_n$ of $M$ with $\Sigma_{n+1} \subset I^+
(\Sigma_n)$ for all $n$ and $\bigcap I^-(\Sigma_n ) = M$ there is
a sequence of $n$-barriers $B_n \subset I^+(\Sigma_n) $ for any $n
\in \N$, then there is a Cauchy surface $\Sigma$ of $M$ such that
$F:= I^+(\Sigma)$ can be foliated by a CMC foliation and the mean
curvature is a temporal function on $F$.
\end{itemize}

In the last item, some results on uniqueness of CMC hypersurfaces appear implicitly, but this question deserves a bigger attention.
A neat problem on uniqueness can be stated as follows, see \cite{ARS}. Consider a Riemannian $n$-manifold $(M,g_R)$, a smooth positive function defined on some interval $f:I\subset \R\rightarrow \R$, and the solutions $u$ to the differential equation on $M$:
\begin{equation}
\begin{array}{c}\label{ebernstein}
u(M)\subset I, \qquad |\nabla u|< f(u) \\
\hbox{div}\left( \frac{\nabla u}{f(u)\sqrt{f(u)^2-|\nabla u|^2}} \right) = nH -\frac{f'(u)}{\sqrt{f(u)^2-|\nabla u|^2}}\left( n + \frac{|\nabla u|^2}{f(u)^2}\right)
\end{array}
\end{equation}
for some constant $H$. The graphs of its solutions can be regarded as the spacelike\footnote{Because of the gradient condition $|\nabla u|< f(u)$.} hypersurfaces  of CMC equal to $H$ in a Generalized Robertson-Walker spacetime $I\times_f M$ (recall, abusing of the notation $g\equiv -dt^2 + f(t)^2g_R$). Notice that this {\em Calabi-Bernstein equation} is the Euler Lagrange one for
the functional $$\mathcal{A}(u)=\int_M f(u)^{n-1}\sqrt{f(u)^2-|\nabla u|^2}dV$$ under the constraint $\int_M\left( \int_{u_0}^uf(t)^n\right)dV=$ constant.
 A specially relevant case of this equation was solved by Cheng and Yau \cite{ChYa} (after the solution by Calabi for  $n\leq 4$):
\begin{theo}
The only entire solutions to Calabi-Bernstein equation \eqref{ebernstein} in $\LL^{n+1}$ 
(i.e., $(M,g_R)\equiv \R^{n}$, $I\equiv \R$, $f\equiv 1$) are linear (or affine) 
 functions.

As a consequence, the only complete maximal hypersurfaces in  
Lorentz-Minkowski space  are the spacelike hyperplanes.
\end{theo}
In fact, they proved that any maximal spacelike hypersurface which is also a closed subset in $\LL^{n+1}$  is a  hyperplane. This yields a surprisingly simple solution to the Calabi-Bernstein problem (recall, for example, that  the analogous results change dramatically for minimal hypersurfaces with dimension larger than seven).

There were, however, well-known counterexamples for the case of CMC hypersurfaces with $H \neq 0$, \cite{St, Tr}. So,
 a line of results about the uniqueness of CMC hypersurfaces has appeared. By using integral inequalities, one can check that   all compact CMC hypersurfaces $\Sigma$ in any GRW spacetime under the {\em null convergence condition} (i.e., Ric$(v,v)=0$ on null vectors) are totally umbilical --- and thus, under mild conditions, $\Sigma$ is a  slice $t=$constant. Such a type of result can be also generalized further, see \cite{Mo, ARS} and references therein.
The feed-back of these results with Riemannian ones have been especially fruitful. Remarkably, both the GRW structure and the restriction for the hypersurfaces of being {\em spacelike}, yields simplifications that inspired some hypotheses for the Riemannian case.

 Under some conditions, previous results can be extended to the case when $\Sigma$ is complete but non-compact, \cite{AlMo}.
In general, for the non-compact case, 
the so-called {\em Omori-Yau maximum principle} (or  {\em asymptotic Cheng-Yau principle}), becomes useful. This principle is stated for complete, connected, noncompact Riemannian manifolds and, roughly speaking, means that any smooth function  $u$ bounded from above on $M$, will admit a sequence $\{x_k\} \subset M$ which plays the role of a maximum (say, $\lim_{k\rightarrow \infty} u(x_k)=$sup$_Mu$,  $|\nabla u(x_k)| \leq 1/k$ and  $\Delta u(x_k)\leq 1/k$).
 The principle holds when the  Ricci curvature is bounded from below as well as in other more refined cases. We refer to  \cite{PRS, AIRtams} and references therein for the recent progress on the Omori-Yau principle and its applications to hypersurfaces in  both, Riemannian and Lorentzian Geometry.

 It is also worth pointing out that, when dim(M)= 3 (i.e., the hypersurface $\Sigma$ is a surface), new tools appear.   For example,  a  different approach to the non-compact case
has been developed very  recently for CMC spacelike surfaces in certain
3-dimensional GRW spacetimes $I\times_fF$; the main idea is to prove that,
 under some natural
assumptions, a metric conformal to the induced one on the surface
$\Sigma$ must be  parabolic, see \cite{RRS} and references
therein. Recall also that the analog to the classical Bj\"orling
problem (construct a minimal surface in $\R^3$ containing a
prescribed analytic strip, solved by H. A. Schwarz in 1890) has
been also considered in the Lorentzian case; this yields a
representation formula for maximal surfaces and allows to
construct new ones explicitly; see \cite{ACM, CDM} for the case of
$\LL^3$ and \cite{oM1} for the general
problem in arbitrary spacetimes, without restriction of the
dimension.

\section{Geodesics and singularity theorems}
\label{Geodesics}

In some concrete spacetimes, singularities might be defined ``by
hand'' but a general definition is difficult \cite{GeWhat}, for example:

\ben

\item The singularity will not be a point of the spacetime, but placed
``at infinity'' -but no natural notion of infinity exists in
general.
 \item The curvature tensor $R$ is expected to diverge, but
all its scalar invariants ($\sum R_{ijkl}R^{ijkl}, \sum \nabla_s
R_{ijkl} \nabla^s R^{ijkl}, S$...) may vanish
 when $R\neq 0$.
\een At any case,  some sort of  ``strange disappearance'' happens
if the spacetime is {\em inextensible, but  an incomplete causal
geodesic exists}, and these two conditions will be regarded as
{\em sufficient} for the existence of a singularity. Then, the aim
of the so-called {\em singularity theorems} is to prove that
causal incompleteness occurs under general natural conditions on
$T$ (an energy condition) and on the causality of the manifold, as
global hyperbolicity. Nevertheless, recall that,  rather than
``singularity'' results, they may be ``incompleteness'' ones: the
physical conclusion of these theorems could be that a physically
realistic spacetime cannot be globally hyperbolic, rather than
being singular.  So, they become ``true singularity'' results when
an assumption as global hyperbolicity is removed... or if SCCC
(Remark \ref{rSCCC}) is true!

Recall the following Hawking's singularity theorem  (see \cite{HE}
or \cite{O} for a detailed exposition):

\begin{theo} \label{thaw} Let $(M,g)$ be a spacetime such that:

1. It is globally hyperbolic.

2. Some spacelike Cauchy hypersurface $\Sigma$ strictly expanding
$H \geq C >0$ ($H$ is the future mean curvature and expansion
means ``on average'')

3. Strong energy (i.e., timelike convergence condition) 
holds: Ric$(v,v)\geq 0$ for timelike $v$.

\noindent Then, any
past-directed timelike geodesic  $\gamma$ is incomplete.
\end{theo}

\noindent {\it Sketch of proof}. 
The last two hypotheses imply that  any past-directed geodesic
$\rho$ normal to $\Sigma$ contains a focal point if it has length
$L'\geq \frac{1}{C}$. Thus, once $\Sigma$ is crossed, no  $\gamma$
can have a point $p$ at length $L>\frac{1}{C}$  (otherwise,  a
length-maximizing timelike geodesic from $p$ to $\Sigma$ with
length $L'\geq L$ would exist by global hyperbolicity, a
contradiction). \qed

\bigskip

This result is very appealing from a physical viewpoint, because
the assumption on expansion seems completely justified by
astronomical observations. Remarkably, the hypothesis $H \geq C
>0$ for some constant $C$ cannot be weakened into $H>0$, as shown
by a surprising (as physically realistic and  far-from-vacuum)
example due to Senovilla \cite{Seno1}. From a mathematical
viewpoint, the reader can appreciate the isomorphic role of the
hypotheses above with the typical ones in Myers type results, say:
global hyperbolicity/ (Riemannian) completenes and timelike
convergence condition/ positive lower bound on the Ricci tensor.

Singularity theorems combine previous ideas with (highly
non-trivial) elements of Causality. Essentially, there are two
types:

\ben \item Proving the existence of an incomplete timelike
geodesics in a global, cosmological setting.

This  is the case of Theorem \ref{thaw}, and some hypotheses there
(specially glob. hyp.) are weakened or replaced by others. For
example, Hawking himself proved that, if $\Sigma$ is compact,
global hyperbolicity can be replaced by assuming that $\Sigma$ is
achronal (i.e., non-crossed twice by a timelike curve). In this
case, the timelike incompleteness conclusion holds, but in a less
strong sense: at least one timelike incomplete geodesic exist.

\item Proving the  existence of an incomplete lightlike geodesic
in the (semilocal) context of gravitational collapse and black holes. \een For
the latter,  the notion of  (closed, future)  {\em trapped}
 surface $K$ (or $n-2$ submanifold)
becomes fundamental. Its mathematically simplest definition says
that $K$ is a compact embedded spacelike surface without boundary,
such that its mean curvature vector field $\stackrel{\rightarrow}{
H}$ is future-directed and timelike on all $K$ \cite{MaSe}
--essentially, this means that the the area of any portion of $K$
is initially decreasing along {\em any} future evolution; when it
is only non-increasing, $K$ will be said {\em weakly trapped}.
Trapped surfaces are implied by spherical gravitational collapse.
One would expect that, at least in asymptotically flat spacetimes
(see next section), they must appear if enough matter is condensed
in a small region
and, under suitable conditions, must imply the existence of a
  {\em black hole} (see \cite{Da} and references therein).
That is, the physical claim is that ``gravitational collapse
implies incompleteness'', and a support for this claim is provided
by the following Penrose's theorem (the first modern singularity
theorem \cite{Pe} -after the works by Raychaudhuri and Komar):

\begin{theo}  Let $(M,g)$ be a spacetime such that:

\ben \item Admits a non-compact Cauchy hypersurface. \item
Contains a trapped surface. \item Ric$(k,k)\geq 0$ for lightlike
$k$. \een Then there exist an incomplete future-directed lightlike
geodesic.
\end{theo}
As emphasized by Senovilla \cite{Seno2}, the  pattern of a
singularity theorem has three ingredients: firstly, a bound on the
Ricci curvature, secondly, a causality condition, and thirdly, an
initial condition on a nonzero-codimensional subset. Remarkably, a
unified treatment of both types of singularity theorems has been
carried out recently by  Galloway and Senovilla  \cite{GaSe}.
Singularity theorems are very accurate, even though it would be
desirable to obtain general results on the nature of the
incompleteness, or ensuring divergences of $R$ in some natural
sense. So, the finding of further types of singularity theorems would be very desirable for physical purposes \cite{Seno1, Seno2}.

\begin{rema} {\em
The subtleties of Lorentzian completeness also appear in the
Lorentzian analogue of Cheeger-Gromoll theorem (see for example
\cite[Chapter 14]{BEE}). To obtain the Lorentzian splitting, of a
spacetime $(M,g)$, dim$(M)>2$, as a product $(\R\times \Sigma,
-dt^2+g_\Sigma)$, where $(\Sigma,g_\Sigma)$ is a complete
Riemannian manifold, one imposes: (a) {\em either geodesic
completeness or global hyperbolicity}, (b) the {\em timelike
convergence condition} (as the meaningful weakening of the
positive semi-definite character of the Ricci tensor in the
Riemannian case),  and (c) the existence of a complete {\em
timelike} geodesic line. }\end{rema}

It is also worth pointing out that the variational approach for
Riemannian geodesics can be extended to the Lorentzian setting but
important particularities appear. It is known since the old work
by Uhlenbeck \cite{Uh} that Morse theory can be applied to
lightlike geodesics, under some conditions (including in particular
global hyperbolicity). Lightlike geodesics satisfy also a
relativistic {\em Fermat principle} \cite{Ko,vP}. Combining both
facts, one can study {\em gravitational lensing}, that is, the
reception at some point $p$ of the spacetime of  light rays
arriving in different directions from the same stellar object, the
latter represented by some timelike curve $c$, see \cite{Pe90,
Pe90bis}. As shown in \cite{CGS}, a very precise result on the
existence and multiplicity of light rays from $c$ to $p$ in
physically realistic spacetimes, can be stated in terms of the
geodesics connecting two points for an appropriate Finsler metric.
For the variational study of spacelike geodesics, see \cite{Ma,
CanSan} and references therein.

\section{Mass, Penrose inequality and CCC}

 {\em Asymptotically flat} 4-spacetimes are useful to model the spacetime around an isolated body.  They can be defined  in terms of Penrose conformal
 embeddings, even though the definition is somewhat involved (see for example \cite{Wa, LRa}). Nevertheless, in what follows it is enough to bear in mind
  that, in an asymptotically  flat (4-)spacetime there exists  spacelike Cauchy
hypersurfaces $\Sigma$
which  admits an {\em asymptotically flat} chart $(\Sigma\backslash K, (x_1,x_2,x_3))$ as follows. 
For  some compact  $K\subset \Sigma$ and some closed ball
$\overline{B_0(R)}$ of $\R^3$, $\Sigma\backslash K$ is isometric
to $\R^3 \backslash \overline{B_0(R)}$ endowed with the metric:

\begin{equation}\label{eh}
h_{ij} - \delta_{ij} \in  O(1/r), \quad \partial_k h_{ij} \in O(1/r^2),
\quad \partial_{k}\partial_l h_{ij} \in  O(1/r^3),\end{equation} in
Cartesian coordinates (this means that $\Sigma$  is intrinsically
asymptotically flat, as a Riemannian 3-manifold; in particular,
Ric and $S$, are in $ O(1/r^3)$), and, even more,
its second fundamental form $\sigma$ satisfies: $\sigma_{ij} \in
O(1/r^2), \partial_k\sigma_{ij} \in O(1/r^3)$.
(This definition can be extended to include more than one end, each one isometric to $(\Sigma\backslash K, (x_1,x_2,x_3))$ as above.)

 The total ADM (Arnowit, Deser, Misner)  {\em mass of an asymptotically flat  Riemannian 3-manifold}  can be defined as 
the limit in any asymptotic chart: 

\begin{equation}\label{emass}  m= \frac{1}{16 \pi} \lim_{r\rightarrow \infty}
\sum_{i,j=1}^3 \int_{S_r} \left(\partial_i h_{ij} - \partial_j
h_{ii}\right) n^j dA, \end{equation} where $n$ is the outward
unit vector to $S_r$, the sphere of radius $r$. Notice that $m$
depends only on the Riemannian 3-manifold; in fact, when this
manifold is seen as a hypersurface of an asymptotically flat
spacetime, the appropriate name for $m$ is {\em ADM energy}, and
the definition of mass depends on $\sigma$, see the next section.
This definition of mass is not mathematically elegant, but recall:

\ben \item  ADM mass appears naturally in a Hamiltonian approach,
as an asymptotic boundary term for the variations of $\int S$.
The definition is not trivial because no strictly local notion of
relativistic energy is available
--- nevertheless,  it is worth pointing out  the attempts to define
a quasilocal mass \cite{LRb}. \item There exists a classical
Newtonian analog when the spacetime is Ricci-flat outside $\R
\times K$, $K$ compact, and there exists a timelike Killing vector
field $\xi$ with lim$_{r\rightarrow
\infty} |\xi| = 1$, such that $\Sigma \perp \xi$. In this case, the divergence theorem yields: 

$m= \frac{1}{4\pi} \int_K
|\xi|^{-1}$Ric$(\xi,\xi) dV = 
\int_K |\xi| \rho dV$

i.e., the
``integral of  the poissonian  density $\rho$ 
measured at $\infty$''. \item The expression in coordinates for
$m$ is manageable: \bit \item  If $h_{ij}=u^4\delta_{ij}$ with
$u(x)= a + \frac{b}{|x|} + O(\frac{1}{|x|^2})$ then $  m=2ab$.

In particular, this is the case if $u$ is ``harmonically flat''
i.e. harmonic with finite limit at $\infty$. \item Otherwise, when
$S\geq 0$ then  $h$ is perturbable to the harmonically flat case
with arbitrarily small error for $m$  and preserving $S\geq 0$
(Schoen and Yau \cite{SY73}; Corvino  \cite{Co} extended the
result for $m>0$ without error in the mass). \eit
\item  Classical outer Schwarzschild metric can be written as: 

$M= \R \times \Sigma$, where $\Sigma=$ $\R^3\backslash
\overline{B_0(|m|/2)}$;

$g= -\left((1- \frac{m}{2|x|})/u\right)^2 dt^2 $ $+$ $h$,
$h_{ij}=u^4\delta_{ij} $ with $u=1+ \frac{m}{2|x|}$

(in particular $\sigma\equiv 0$). Of course, the classical
Schwarzschild mass $m$ agrees ADM mass. \een One expects from the
physical background that, when the dominant property holds, the
ADM mass will be positive for any asymptotically flat Cauchy
$\Sigma$.
Two technical points are relevant here: (a)  When $\Sigma$ is
totally geodesic ($\sigma\equiv 0$)  the dominant property
yields 
$S\geq 0$. (b) Under our definition of asymptotic flatness,
$\Sigma$ is necessarily complete, but  the Riemannian part of
exterior Schwarzschild spacetime
$(\R^3\backslash\overline{B_0(|m|/2)},h)$ is incomplete for any
$m\neq 0$. Of course, this is not a problem for the computation of
the limit in the expression of the ADM mass, and one can also
extend and modify $(\R^3\backslash\overline{B_0(|m|/2)},h)$ in a
bounded region to obtain a complete Riemannian manifold $\Sigma^c$
with the same asymptotic behaviour. Moreover, in the globally
hyperbolic case $m > 0$, one can obtain such a $\Sigma^c$ (say,
corresponding to the spacetime created by a star of the same mass)
with: (i) the same asymptotic behaviour, (ii) $S\geq 0$. Clearly,
this property is not expected in the non-globally hyperbolic case
$m<0$. And, in fact, it is forbidden by
the  {\em Riemann positive mass theorem}:

\begin{theo}  Let $(\Sigma, h)$ be
any asymptotically flat (complete) Riemannian manifold with $S\geq
0$. Then, $m\geq 0$ and equality holds iff $(\Sigma, h)$ is
Euclidean space $\E^3=(\R^3,\delta).$
\end{theo}
\begin{rema} {\em This celebrated result by Schoen and Yau \cite{SY72}
is a purely Riemannian one. From this case, more general
``positive mass'' results follow, which include the case $\sigma
\not\equiv 0$ \cite{SY73}; see also the comments in the next
section about Witten's, completely different, proof.
 By the way, recall that the solution of Yamabe problem was completed by
using the above result (see the nice survey \cite{LP}).

It is worth pointing out that, because of a technical problem
which goes back to the known failure of regularity of minimal
surfaces in dimensions greater than 7, the positive mass theorem
is proved for dimensions up to 7 (this has been completed only
recently, by using Schoen and Yau techniques, see \cite{EHLS}) as
well as for spin manifolds of any dimension (by using Witten's
techniques to be explained in the next section). }\end{rema}

Next, we will consider a no less spectacular further step (for a
detailed exposition, see \cite{Br2}). But, first two notions will
be briefly explained:

1.- {\em WCCC}. A question related with SCCC (see Remark
\ref{rSCCC}) is the so-called {\em weak cosmic censorship
conjecture } (WCCC), which is stated in the framework of
asymptotically flat spacetimes. In such spacetimes, a natural
notion of asymptotic future null infinity ${\cal J^+}$ can be
defined (${\cal J^+}$ is a subset  of the image of $M$ for a
suitable conformal embedding in a bigger spacetime $\bar M$) and,
then, also a rigorous notion of the {\em black hole} region $B$ of
$M$ appears ($B= M\backslash J^-({\cal J^+}$)) ---this region
corresponds to the intuitive idea of a ``spatially bounded region
from where nothing can escape''. WCCC asserts that (maybe only
generically)  any spacetime $M$ obtained as the maximal evolution
of physically reasonable initial data with an asymptotic
decay\footnote{Typically, this data must satisfy: (i) $(\Sigma, h,
\sigma)$ is asymptotically flat, (ii) $T$  satisfies the dominant
property, and the equations for $T$ constitute a quasilinear,
diagonal, second order hyperbolic system, (iii)  the fall-off of
the initial value of $T$ on $\Sigma$ is fast enough for the
$h$-distance, and $h$ is assumed to be complete.},
will be asymptotically flat and, in a restrictive sense, globally
hyperbolic at infinity\footnote{More precisely, the latter means
that the spacetime is {\em strongly asymptotically predictable},
see \cite{Wa}. Recall that WCCC cannot be regarded as a particular
case of SCCC.}.
 The physical interpretation of this assertion is that no
singularity (except at most an ``initial'' one) can be observed
from $M\backslash B$, that is,
 singularities must lie inside a black hole and cannot be seen from outside
 (singularities are not ``naked'').

2.- {\em Outermost trapped surfaces}. 
Given a totally geodesic asymptotically  flat slice $\Sigma$,
those trapped surfaces (more precisely, compact spacelike surfaces
whose expansion respect to the outer future lightlike direction is
at no point positive) contained in $\Sigma$  which are boundaries
of a 3-manifold, are known to satisfy:
\ben \item Such trapped surfaces correspond to compact minimal
surfaces of $\Sigma$.

\item The  outermost boundary compact minimal surfaces
(necessarily  topological spheres, each one the ``apparent horizon in $\Sigma$ of a black hole'') are well-defined. 
\item Let ${\cal H}$ be the union of the outermost minimal
surfaces. Under WCCC, if ${\cal H}$ is connected and $A_0$ denotes
its area, physical considerations ensure that the ``contribution
to the mass'' $m_0$ of the corresponding black hole would satisfy:
$m_0 \geq \sqrt{\frac{A_0}{16 \pi}}$.
\een 
Therefore, choosing any asymptotic $\Sigma$ one expects for its
mass $m_\Sigma$: \beq
 \label{star} m_\Sigma \geq
\sqrt{\frac{A_0}{16 \pi}}\quad \quad \eeq
 (at least if the second fundamental form vanishes).
But (\ref{star}) is an inequality in  pure Riemannian Geometry.
Thus, the following precise result must hold:

\begin{theo}
Let $(\Sigma, h)$ be an asymptotically flat Riemannian 3-manifold
with $S\geq
0$, and 
 let ${\cal
H}_0$ be the largest outermost (connected)  minimal surface, 
with area $A_0$. Then inequality (\ref{star}) is satisfied,
and the equality holds if and only if $(\Sigma, h)$ is
Schwarzschild Riemannian metric outside ${\cal H}_0$.
\end{theo}
This is the celebrated ``Riemann-Penrose inequality'', proved by
Huisken and Ilmanen \cite{HI} (who re-prove then the Riemann
positive mass theorem), and shortly after extended by Bray to the
full area  of the (maybe non-connected)
${\cal H}$,
with a different proof \cite{Br} based on positive mass theorem.

Penrose inequality is a more general conjecture, which includes
the full spacetime case $\sigma \neq 0$ (recall that  the case
above would correspond when $(\Sigma, h, \sigma = 0)$ can be
regarded as an initial data set for the spacetime).
It is still open, and it becomes  a major problem in Differential
Geometry. An evidence of its difficulty is that  it is supported
by physical grounds, and counterexamples to more general appealing
mathematical conjectures have been found, see \cite{CaMa}; we
refer to the reviews \cite{CGP, Mars} for comprehensive
references.

\section{Spinors and holonomy}
\label{Spinors}

Dirac operators are popular objects of study in the area of global analysis, one of the main reasons being the existence of index theorems for them, see the standard textbook by Lawson and Michelson \cite{LM89} or the book by Berline, Getzler, Vergne \cite{BGV} on Dirac operators, both almost exclusively treating the Riemannian situation --- this reflects the fact that index theory presently is applicable almost exclusively to elliptic operators. Another main reason of interest in spinors, more predominant in the Lorentzian case, is the presence of {\em Weitzenb\"ock formulas}. These formulas reflect the fact that the Dirac operator as a natural first-order operator on spinors is a root of the Laplacian type operator plus a curvature-induced zeroth order term.

While in the book by Lawson-Michelson real spinors play a prominent role, in the following, we want to focus here on {\em complex} spinors.

Spinor bundles are defined verbatim in the same way as in the Riemannian case, with $SO(n)$ always replaced by the connected component of the identity of $SO(1,n)$, but there are some important differences of the Lorentzian to the Riemannian case: The natural (pseudo-Hermitean) scalar product $\langle \cdot , \cdot \rangle $ on the spinor bundle is not definite, but of split signature. Any timelike vector field $X$ can be used to define a (non-natural) positive-definite scalar product $( \cdot, \cdot ) := \langle X \cdot  , \cdot ) $. Clifford multiplication is $\langle \rangle$-symmetric (instead of antisymmetric, as in the Riemannian case).

For a pseudo-Riemannian spin manifold of arbitrary signature, one
can define the {\em Dirac operator} on $C^1$ (or at least $W^{1,p}$) sections $\psi$ of the spinor bundle
by $D \psi := \sum_{i =1}^n \epsilon_i e_i \cdot \nabla_{e_i} \psi $ (where the $e_i$ are a
pseudoorthogonal basis and $\epsilon_i$ is the sign of $g(e_i,
e_i)$). The Dirac operator is formally self-adjoint, essentially
self-adjoint if $(M,g)$ is complete, and satisfies the
Weitzenb\"ock identity

$$   D^2 = \nabla^* \nabla + \frac{1}{4} S  , $$

\v where $S$ is the pseudo-Riemannian scalar curvature of $(M,g)$.

In the Riemannian situation, as the connection Laplacian $\nabla^*
\nabla $ is positive-definite, the Weitzenb\"ock formula is the
initial point of many obstructions to positive scalar curvature
for spin manifolds. The Weitzenb\"ock formula for the Lorentzian
Dirac operator looks superficially the same but the connection
Laplacian here is a {\em hyperbolic} operator instead of an
elliptic operator.

Exactly as in the Riemannian situation, any spinor defines an associated one-form, the so-called {\em Dirac current}. In Lorentzian geometry, an additional factor $i$ appears in the definition, basically to balance the effect of the aforementioned differences. As in the Riemannian case, elementary calculations show that the Dirac current of a parallel spinor is a parallel spinor field. While in Riemannian geometry, the Dirac current of a {\em real} Killing spinor is a Killing vector field, in the Lorentzian case the same is true for {\em imaginary} Killing spinors. In stark contrast to the Riemannian situation, in Lorentzian geometry the Dirac current is always non-trivial for a non-vanishing vector field. The Dirac current of any eigenspinor of a twisted Dirac operator is always divergence-free, and $(\cdot, \cdot)$ can be used to define a conserved charge.

Many properties of spacetimes carrying special spinor fields can be read off from their Dirac current. E.g., as shown by Ehlers and Kundt \cite{EK62}, a four-dimensional Lorentzian spin manifold with a parallel spinor is locally isometric to a pp-wave (in this case, the Dirac current is a parellel null vector field).

An important application of spinors in Lorentzian geometry is the
commented proof of the positive mass theorem due to the seminal
ideas by Witten \cite{Wi}, made rigorous by Parker and Taubes
\cite{PT} and others (see the independent work by Reula
\cite{oR82} as well as \cite{He} and references therein). In fact,
the spacetime viewpoint is necessary here.  So, we will revisit
the approach in the previous section, and focus in the
 positiveness of the energy (which, as commented above, could be also proven
 by using Schoen and Yau techniques). The starting point is a Cauchy
hypersurface $(\Sigma,g \vert_\Sigma)$ that is asymptotically flat
in the sense explained around formula \eqref{eh}, including the
bounds for the second fundamental form. In this case, the
expression of $m$ in \eqref{emass} is taken as the definition of
the energy $E$,  and the momenta   $P_l$ are defined by:


$$   P_i = \frac{1}{8\pi} {\rm lim}_{r \rightarrow \infty} \int_{S_r}   \big(  \sum_{j=1}^3 \sigma_{ij} n^j - \sum_{j=1}^3 \sigma_j^j  n_i \big)  dA , $$

\v where $\sum_{j=1}^3 \sigma_j^j/3$ is the mean curvature of
$\Sigma$. As we have seen, $E$ is independent of the chosen
asymptotic coordinates, and the freedom in the choice of these
coordinates yields new momenta $P'_1, P'_2, P'_3$ which differ in
an element of $O(3)$. This allows to construct a vector $ V:= (E,
P_1, P_2, P_3) \in \R^{1,3}$, the {\em  ADM energy-momentum} (a
different choice $\Sigma '$ of hypersurface would yield a new
vector $V'$ which would be related to $V$).  The statement of the
positive energy theorem is that $V$ is a causal vector and equal
to $0$ if and only if the spacetime is flat around $\Sigma$.
Witten's idea was to use a Weitzenb\"ock formula for the spacetime
Dirac operator applied to spinors tangent to the hypersurface and
extended parallely along normal geodesics in a small normal
neighborhood of the Cauchy hypersurface. The dominant energy
condition then ensures that the zeroth-order term in the
Weitzenb\"ock formula is positive, that yields directly the
positiveness of the energy. Moreover, one
 can find a harmonic spinor approaching in coordinates a
parallel spinor on $\R^{3}$  (thought of as in embedded via the
Cauchy hypersurface) such that the limit of the boundary term
appearing in the integral form of the Weitzenb\"ock formula is
exactly $ E - \vert P \vert$, thus obtaining $E\geq |P|$ i.e., the
energy momentum is a causal vector (a corollary is then the
positive mass theorem $E \geq 0$). A central tool to prove
existence of these solutions are Green's functions in weighted
Sobolev spaces performed in detail by Parker and Taubes. In 1987,
Yip \cite{pY87} showed that the energy-momentum vector has to be
even {\em timelike} (non-lightlike) unless $M$ is flat around the
Cauchy hypersurface, also by using spinors techniques. As already
pointed out, Eichmair et al.  \cite{EHLS} have given a recent
proof of the energy-momentum inequality $E\geq |P|$ in the case
that the manifold $M$ is not necessarily spin, but that ${\rm dim}
(M) \leq 7$. This is obtained by using Schoen and Yau techniques
but, remarkably, a new difficulty appears, as minimal surfaces are
now replaced by marginally outer trapped hypersurfaces, which do
not come from a variational characterization.



Another fundamental concept in geometry connected to spinors is
{\em holonomy}. That notion can be defined on any bundle with a
connection, denoting, at a point $p$, the group of diffeomorphisms
of the fiber over $p$ which are parallel transports along curves
starting and ending at $p$. If applied to a semi-Riemannian
manifold of signature $(m,n)$ with its Levi- Civita connection,
it is a restriction of the standard representation of $SO(m,n)$ to a subgroup.
It is easy to see that in a connected manifold, the equivalence class of the holonomy representation does
not depend on the point $p$. The corresponding infinitesimal
notion (taking the Lie algebra of the holonomy group) is called
the {\em holonomy algebra}.

Now, in the Riemannian case, we have Berger's list: a
simply-connected Riemannian manifold is either locally symmetric
or can be decomposed as a Riemannian product each $k$-dimensional
factor of which has the holonomy $SO(k) $, $U(k/2)$, $SU(k/2)$,
$Sp(k/4) \cdot SP(1)$, $Sp(k/4)$, $G_2$ (in which case $k=7$) or
$Spin(7)$ (in which case $k=8$). The Lorentzian case is a bit more
involved and has remained open until recently. One of the
difficulties compared to the Riemannian case is the difference
between decomposability and reducibility pointed out in
Section \ref{s2}. In fact, classical de Rham Riemannian
decomposition relies on the fact that the orthogonal complement
$A^\perp$ of an invariant subspace $A$ must be not only invariant
too, but also a complement of $A$. When the latter property holds
in the semi-Riemannian case, an analogous de Rham-Wu decomposition
is obtained, but this is not the case  when $A$ is degenerate
--- something that can occur in the Lorentzian case. So, the
elementary building blocks of the Riemannian classification,
irreducible subspaces, have therefore to be complemented by new,
properly Lorentzian, building blocks, the indecomposable but
non-irreducible subspaces. Such an $m$-dimensional subspace
contains an invariant lightlike subspace $N$, and its holonomy
algebra is contained in $(\R \oplus so(m-2) ) \ltimes \R^{m-2}$,
thus a central part of the classification is done by the
$so(m)$-projections of the possible holonomy representation, the
so-called {\em screen holonomy} acting on the associated invariant
codimension-2 subbundle of $\tau_M$, the {\em screen bundle} given
by $N^{\perp}/N$ which carries a well-defined Riemannian metric.
Thomas Leister in his PhD thesis in 2005 (published in
\cite{WASISTDAS}) showed that the screen holonomy is always the
holonomy algebra of a Riemannian manifold (and, as remarked by
Anton Galaev, this is an exceptional feature of Lorentzian
geometry not present in higher signatures). In that manner, he was
able to solve some of the remaining problems in the Lorentzian
classification by means of the corresponding Riemannian
techniques, and finally obtained the full classification. Galaev
\cite{FUNKTIONIERE} gave then analytic examples for every holonomy
representation of Leistner's list. Still, a missing piece were
examples of {\em globally hyperbolic} manifolds {\em with complete
Cauchy hypersurfaces} with the given holonomy representations.
This was done for manifolds with parallel spinors (for which case
the above classification yields groups $G \semidir \R^n$ for $G$
being a product of $SU(p), Sp(q), G_2$ or $Spin(7)$) in an article
of Helga Baum and Olaf M\"uller \cite{hBoM06}, via a cylinder
construction analogous to one by B\"ar, Gauduchon and Moroianu
\cite{BGM}, building a parallel spinor from a so-called Codazzi
spinor, and another construction relating Codazzi spinors to
imaginary Killing spinors (whose importance in geometry is
explained in the next paragraph). In 2013, Helga Baum and Thomas
Leistner then solved the analytic initial value problem for
parallel spinors \cite{BL>13}.


A {\em Killing spinor} is a  spinor $\psi$ such that there is a
constant $b \in \C$ with $ \nabla_X \psi = b X \cdot \psi$ for all
vectors $X$. As shown by Friedrich in \cite{tF}, Killing spinors
can serve as landmarks where spectral estimates get sharp, in the
following sense: If $M$ is compact and the scalar curvature is
bounded from below by a positive constant $s_0$, then for all
eigenvalues $a$ of the Dirac operator we have $a^2 \geq
\frac{1}{4} \frac{n}{n-1} s_0$, and equality in this estimate
implies that the corresponding eigenspinor is a Killing spinor.
One can consider the modified connection $\nabla_b := \nabla - b
\1 $ to conclude that a Killing spinor never vanishes. An
elementary calculation shows that $($Ric$(X) - 4 b^2 (n-1) X)
\cdot \psi =0$, and that implies (by nonvanishing of $\psi$) that
the image Ric$- 4 b^2 (n-1) \1_{TM}$ is contained in the null
cone. Taking the trace once more, one sees that the scalar
curvature equals $4n (n-1) b^2$, in particular, $b$ is either real
or purely imaginary. A Killing spinor is called real resp.
imaginary if $b \in \R $ resp. $b \in i \R$. Over several years,
different people aimed at a full classification of Killing
spinors. Christian B\"ar \cite{cB} finally came up with a cone
construction which associated to each Killing spinor on a manifold
$M$ a parallel spinor on the Riemannian cone over $M$. As the
existence of parallel spinors leads to special holonomy, B\"ar
obtained a classification via the classification of Riemannian
holonomies. Imaginary Killing spinors were classified by Helga
Baum in \cite{hB89-1}, \cite{hB89-2} in a completely different
way: Let $(M,g)$ be a complete connected spin manifold. It carries
an $i \cdot a$-Killing spinor iff it is a warped product $\R
\times_{e^{-4at}} N $ for a complete connected spin manifold with
a non-zero parallel spinor field. The idea is to show that the
manifold is foliated by level sets of the norm $t$ of the spinor
field. Christoph Bohle \cite{cBo} examined real Killing spinors on
Lorentzian manifolds, relating them also to warped products.
Felipe Leitner \cite{fL}, finally, considered imaginary Killing
spinors on Lorentzian manifolds. Their Dirac current is easily
seen to be causal, and when it is null, then the manifold is
Einstein.

\section{Some further topics and a double invitation}

In this article, as announced, we intend
to invite experts from other branches of mathematics, especially
in Riemannian geometry and global analysis, in two respects: Firstly, we invite {\em
users} of Lorentzian geometry. We hope that the article made clear that
Lorentzian geometry can be extremely useful not only in physics,
but also in mathematical contexts. One famous example is the
aforementioned solution of the Riemannian Yamabe problem via
Lorentzian techniques. To a large extent, this potential of
Lorentzian geometry remains unexplored up to now.

\bigskip

Secondly, we want to invite {\em providers}. The open topics
in Lorentzian geometry do need support from other branches of
mathematics. In the following, we list some important open
questions in Global Lorentzian Geometry (without any claim of
completeness). In order to do so, it is convenient to distinguish
between those that arise directly in Mathematical Relativity and
those that are mathematically natural, independent of physical
motivations. Along this paper we have emphasized  some of the
first type. But recall that the questions on Lorentzian manifolds
inspired only in reasons of mathematical naturalness and beauty,
are interesting in their own right and, sooner or later, will have
applications to General Relativity or other parts of Mathematical
Physics --- recall that General Relativity is one of the two
fundamental physical theories, and Quantum Theory the other one. For decades, 
the attempt to unify both theories has been a physical challenge
 and a permanent source of mathematical inspiration.

Along this article some  open questions in Mathematical Relativity
has appeared more or less explicitly, such as: (a)  Cosmic
Censorship Conjecture (weak and strong), including full Penrose
inequality, (b) Cauchy problem (blow up criteria, global
regularity for large data...), or (c) definitively satisfactory
definition of singularities, including both singularity theorems
(which involve divergences on curvature and not merely
incompleteness) and a precise description of the {\em boundary} of
the spacetime. Of course, there are many other relevant questions in
Mathematical Relativity (see \cite{CGP}). We would like to point
out here the interest attracted by the questions of stability.
Christodoulou and Klainerman \cite{CK} proved the non-linear
stability of Lorentz-Minkowski spacetime $\LL^4$ as a solution of
Einstein equation. This means that a small perturbation of the
initial conditions for $\LL^4$ yields a spacetime with properties
close to $\LL^4$ (and, for example, not to a spacetime with
singularities). In spite of the  simplicity of this idea, the
proof is extremely difficult -recall that \cite{CK} is a 500 pages
book. The result is  a landmark in Mathematical Relativity, and
opens the study of the stability under weaker falloff hypotheses
of the initial data or the stability of other spacetimes, as those
with constant curvature, or of the Einstein equation coupled to
other field theories.

Finally, let us point out some more purely mathematical questions,
some of them suggested above, but only tangentially.
(1) Classification of submanifolds 
with natural geometric properties (constant mean curvature,
umbilic, etc.) in spaceforms and other physical or mathematically
relevant spacetimes; notice that some of these questions had
motivations from the viewpoint of the initial value problem and
were commented in Section 5, but such problems evolve further,
independent of physical motivations. (2) Critical curves for
indefinite functionals on Lorentzian manifolds: even though the
role of geodesics in General Relativity gives a general support
for this, the infinite-dimensional variational mathematical
approach for geodesics, including spacelike ones, has independent
interest, see the seminal works by Benci, Fortunato and Giannoni
\cite{BFG}, the book \cite{Ma} of the review \cite{CanSan}; we emphasize
that even a simple question as if any compact Lorentzian
manifold must admit a closed geodesic remains open. (3)
Curvature:
 curvature bounds 
groups   have been stressed above, but there are many other
questions
   related with curvature operators, e.g., those starting at the Osserman problem,
   solved a decade ago, see \cite{Rio}. (4)
Classification of Lorentzian spaceforms: such a topic has a deep
importance and  tradition in Geometry, we
recommend the recent revision of a paper by Mess in \cite{Mess} as
an example of this exciting problem. (5) Links between Lorentzian
and Finslerian geometries at different levels are being developed
fast in the last years, see \cite{CJM, CJS, FHSmemo} as a sampler.

The second-named author is partially supported by the Grants MTM2010--18099 (MICINN) and P09-FQM-4496 (J. Andaluc\'{\i}a) with FEDER funds.

\vspace*{-1.0ex}
\begin{figure}[h]
\begin{minipage}{0.54\textwidth}
\centering
\includegraphics[scale=0.2]{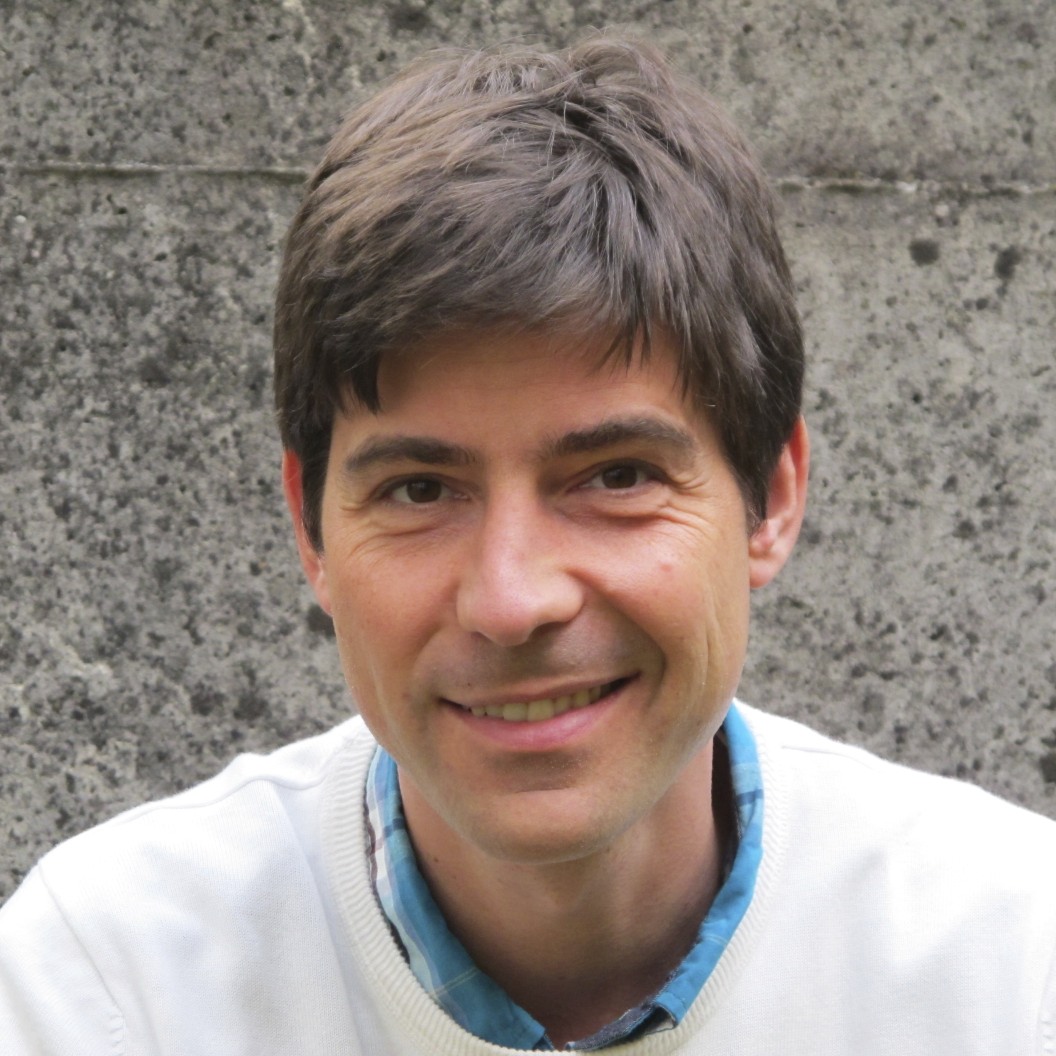}
\end{minipage}
\begin{minipage}{0.45\textwidth} \setlength\parindent{1em}
\centering
\includegraphics[scale=0.3]{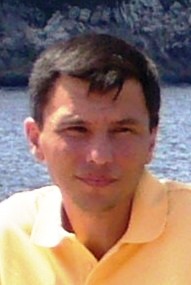}
\end{minipage}
\end{figure}

\begin{minipage}{0.45\textwidth}
Olaf M\"uller obtained his PhD at the Max-Planck Institute for Mathematics in the Sciences, Leipzig, in 2004. Since then, he has worked at the Humboldt University Berlin in the group of Helga Baum, at the Universidad Nacional Aut\'onoma de M\'exico (UNAM) and at the University of Regensburg in the group of Felix Finster. His research focuses on differential geometry and global analysis, in particular the geometry of globally hyperbolic manifolds. 
\end{minipage}
\hfill
\begin{minipage}{0.47\textwidth}
{\footnotesize
Miguel S\'anchez is a mathematician and physicist at the Department of Geometry and Topology of the University of Granada. He obtained his PhD in Mathematics in 1994 and, since then, has worked in several topics on Differential Geometry, Global Analysis on Manifolds and Mathematical Physics. He has supervised four PhD thesis, and his publications include more than fifty articles of research, as well as a couple of text books, one of them about Lorentzian Geometry.}
\end{minipage}

  \end{document}